\documentclass[12pt,a4paper]{article}
\usepackage{amsmath}
\usepackage{amsthm}
\usepackage{amsfonts}
\usepackage{amssymb}
\usepackage{calc}
\usepackage{multicol}
\usepackage{makeidx}
\usepackage{eucal}
\usepackage{exscale}
\usepackage{epic}
\usepackage{eepic}
\usepackage{pstricks}
\usepackage{wrapfig}
\usepackage{color}
\usepackage{rotating}
\usepackage{ecltree}
\begin{document}

\title {Krein's theory on strings applied to fluctuations of L\'evy processes}
\author{Sonia Fourati}
\maketitle
\centerline{\it LMI de l'INSA de Rouen, place Emile Blondel, 76130 Mont St Aignan. France}
\centerline{\it LPMA des Universit\'es Paris VI et VII, 4 Place Jussieu, case 188, 75252 Paris, Cedex 05. France}

\bigskip

{\bf ABSTRACT}
We give an interpretation of the bilateral exit problem for L\'evy processes via the study of an elementary Markov chain. 
We exhibit a strong connection between this problem and Krein's theory on strings. For instance, for symmetric L\'evy processes with bounded variations, the L\'evy exponent is the correspondant spectral density  and the Wiener-Hopf factorization  turns out to be a version of Krein's entropy formula.

\section{Wiener-Hopf factorization conditionned by the amplitude}
Let $(\Omega,{\cal F})$ be the canonical path space of real valued functions with finite or infinite  life time $\zeta$, $X$ the canonical process, $\check X$ the dual process $-X$. Denote by ${\bf P}$ the law of a L\'evy process and ${\bf Q}$ the``law" of   $X$  killed either at an exponential time or at a ``Lebesgue''  time. In other words ${\bf Q}$ is either : 

(``First case") : $${\bf Q}(dw)=\int_0^{+\infty}{\bf P}((X_s)_{s\in[0,t[}\in dw)e^{-t}dt$$ or (``Second case") :  $${\bf Q}(dw)=\int_0^{+\infty} {\bf P}((X_s)_{s\in[0,t[}\in dw)dt$$ 

Let us define $$F:=X^-_{\zeta}\hbox{( = the final value of the canonical process)}$$
$$M:=\sup\{X_t ; t\in[0,\zeta[\}\quad\hbox{ (= the maximum value)}$$
$$m:=\inf\{X_t ; t\in[0,\zeta[\}\quad \hbox{(=the minimum value)}$$
Then $$M-m \hbox{ is the amplitude.}$$

In the sequel, we exclude  compound Poisson processes so that the time $t$ when $X_t$ 
 is at its mimimum  value is unique and denoted by $\rho$.

Let   ${\bf Q}^{\uparrow}$ be the ``law after the minimum ", more precisely,  in the first case 
$${\bf Q}^{\uparrow}(dw)={\bf Q}( (X_{s+\rho}-m)_{s\in [0,\zeta-\rho[}\in dw)$$ 
and  in the second case, 
 $${\bf Q}^{\uparrow}(dw)={\bf P}(\int_0^{\zeta}1_{(X_t-X_{t-s})_{s\in ]0,t]}\in dw}L(dt))$$ 
where  $L$  is a local time at $0$ of the reflected process $\sup_{s\in [0,t]}X_s-X_t$, (${\bf Q}^{\uparrow}(dw)$ is thus defined  up to a multiplicative constant ): 

The following independence property (in both cases)  is well known :
$${\bf Q}( (X^-_{\rho -s}-m)_{s\in [0,\rho[}\in dw_1), (X_{s+\rho}-m)_{s\in [0,\zeta-\rho[}\in dw_2)=\check{{\bf Q}}^{\uparrow}(dw_1){{\bf Q}}^{\uparrow}(dw_2)$$

 \medskip
 
{\bf Proposition 1} {\it In both cases and for every $x>0$ :} 
$${\bf Q}( -m\in du, M-m\leq x, F-m\in dv)=\check{{\bf Q}}^{\uparrow}(M\leq x, F\in du){{\bf Q}}^{\uparrow}(M\leq x, F\in dv)$$

\medskip 

{\bf Argument :} The  event $M-m\leq x$ is the intersection of  the two events  : ``The amplitude before the minimum does not exceed $x$" and ``the amplitude after the minimum does not  exceed $x$". The result then follows easily from the indepandence property given above.
\medskip 

The purpose of the sequel is to give information  on the two laws ${{\bf Q}}^{\uparrow}(M\leq x, F\in dv)$ and $\check{{\bf Q}}^{\uparrow}(M\leq x, F\in du)$ for every $x>0$. They are obviously a caracterization of the bilateral problem or equivalently of the law of the trivariate variable $(-m,M,F)$.

\section{An elementary Markov Chain} 
Denote :
$$T_0=\rho$$ 
$$Z_1=\sup\{X_t-m\ ;\ t > \rho\}\quad \{T_1\}= \{t\  ;\  X_t\vee X_t^- -m=Z_1\}$$
$$Z_2= \inf \{ X_t-X_{T_1}\vee  X_{T_1}^- \ ;\  t > T_1\} \quad \{T_2\}=\{t > T_1\ ;\  X_t\wedge X_t^- -X_{T_1}\vee  X^-_{T_1}=Z_2\}$$
etc. Note that
$$ Z_1>0, Z_2<0,\dots, \qquad |Z_n| \geq |Z_{n+1}|\qquad \sum_{n=1}^{\infty} Z_i=F-m$$

 let $$H(dx):={\bf Q}^{\uparrow}(M\in dx)={\bf Q}( Z_1\in dx)$$
$$H(x):={\bf Q}^{\uparrow}(M\leq x)={\bf Q}(Z_1\leq  x)=\phi(x,0)$$

We obtain easily the next two propositions 
\medskip

{\bf Proposition 2}  :  {\it $(Z_n)$ is a Markov chain whose  transition kernel is }

$$P(x,dy)={{\check H}(-dy)\over {\check  H} (x)}1_{y\in [-x,0]} \quad \hbox{ for } x\geq 0$$
$$P(x,dy)={H(-dy)\over H(x)}1_{y\in [0,-x]} \quad \hbox{ for } x\leq  0$$
\medskip 
Put 
$$\phi(x,\lambda ):={\bf Q}^{\uparrow}(M\leq x,e^{-\lambda F})={\bf Q}(Z_1\leq x, e^{-\lambda \sum_{n=1}^{+\infty}Z_n})$$

\medskip

{\bf Proposition 3 } \ : \  {\it For every $\lambda  \in C$  ,

$$\phi (dx,\lambda)=e^{-\lambda x}{\check \phi}(x,-\lambda){H(dx)\over {\check H}(x)}$$

$${\check \phi }(dx,-\lambda)=e^{\lambda x}{\phi}(x,\lambda){{\check H}(dx)\over H(x)}$$

For $x\to 0^+$ : $$\phi(x,\lambda)\sim H(x) \quad {\check \phi}(x,\lambda)\sim {\check H}(x)$$

For $x\to +\infty$ and $\lambda \in iR$ :
 $$\phi(x,\lambda)={\bf Q}^{\uparrow}(e^{-\lambda F})
\quad \check {\phi}(x,-\lambda)=\check{\bf Q}^{\uparrow}(e^{\lambda F})$$
}
\section {An example : Stable processes killed at a Lebesgue time}

In this section $X$ is a standard stable process under $\bf P$ and 
$${\bf Q}(dw) :=\int_0^{+\infty}{\bf P}((X_s)_{s\in[0,t[}\in dw)dt$$ 
Let $\alpha$ denote the stability index, 
$$\gamma:=\alpha{\bf P}(X_1>0)\qquad \delta :=\alpha{\bf P}(X_1<0)$$
 $$\beta_{a,b}(t) := {\Gamma(a+b)\over \Gamma(a)\Gamma(b)}t^{a-1}(1-t)^{b-1} 1_{t\in [0,1]}\quad {\hat \beta}_{a,b}(\lambda ) :=\int_0^{+\infty}e^{-\lambda  t} \beta_{a,b}(t) dt$$

\medskip 
{\bf Theorem :} {\it For some  constant $k\in ]0,+\infty[$, we have the identity}

$$\phi(x,\lambda )=k.x^{\gamma} {\hat \beta}_{\gamma,\delta+1}(\lambda)$$
\medskip 

{\bf Argument} : It is clear that in the case involved, the functions $H$ and $\check H$ are 
of the form $H(x)=c x^{\gamma}$ and $\check H(x) =\check c x^{\delta}$ ($ c,\check c \in ]0,+\infty[$). When iterating the differential system given in proposition 3, one gets that the function $x\mapsto \phi(x,\lambda )$ satisfies  the confluent hypergeometric equation whose solutions are known.

{\bf Consequences  }(essentially [Rogozin])

1) $${{\bf P}\bigl(\int_0^{+\infty} 1_{X_s-\inf_{s\in [0,t]} X_s\leq 1}dL_s\bigr) \over {\bf P}(\int_0^{+\infty} 1_{X_s\leq 1}dL_s) }={\Gamma(\delta+1)\over \Gamma (\alpha+1)}$$

where $L$denote a  local time at $0$ of the reflected process $\sup_{s\in [0,t]} X_s-X_t \Bigr)$

2) Under ${\bf Q}^{\uparrow}$,   ${F\over M}$ and $M$ are independant and ${F\over M}\stackrel {d}=\beta_{\gamma+1,\delta}(t)dt\hfill$

3) $$\int_0^{+\infty}{\bf P}(S_t\leq x,-I_t\leq y) dt ={\Gamma(\gamma+1)\Gamma(\delta+1)\over \Gamma^2(\alpha
+1)}x^{\gamma}y^{\delta}$$

4) $${\bf P} (T^a<T_b)= \int _0^{b\over a+b} \beta_{\gamma,\delta} (t)dt$$

($T^a=\inf \{ t ;\ X_t>a\}\quad T_b=\inf\{t ;\ X_t<-b\}$).

\section{A survey of Krein theory on strings} 

Put $${\cal M}:=\{m :\ [0,+\infty[\to [0,+\infty[,\  m\nearrow,\ m([0,+\infty[)=[0,+\infty[\}$$

An element of ${\cal M}$ is called a ``string".

$${\cal D}:=\{\Delta \hbox{ measure on }\  [0,+\infty[ \hbox{ such that }$$
$$ \int _1^{+\infty} {\Delta(du)\over u^2}<+\infty; \int _0^1 {\Delta(du)\over u^2}=+\infty ;\  \int_1^{+\infty} \Delta(du)=+\infty\}$$
An element of ${\cal D}$ is called a ``spectral measure". 
 \medskip
 
 \centerline {\bf FACTS}
 \medskip 
 
{\it For every string $m\in {\cal M}$, we have the following properties, 
 
1)  For every $\lambda > 0$, the equation $${d^2\over dm dx}X=\lambda^2 X$$ admits  a  unique solution (called ``A-solution'') with conditions $$A(0)=1\hbox{ and } A'(0)=0$$
and a unique solution (called ``D-solution'') with  conditions $$D'(0)=-1\hbox{ and  }D(+\infty)=0$$
 
 2) $$D(x)=A(x)\int_x^{+\infty}{1\over A^2(t)}dt$$

3) There exists a unique $\Delta \in {\cal D}$, (``the spectral measure of  the string $m$''), such  that  for every $\lambda>0$ :
  
$$D(0,\lambda)={2\over\pi} \int _0^{+\infty} {\Delta(du)\over u^2+\lambda^2}$$

moreover $m\rightarrow \Delta$ is a bijection from ${\cal M}$ onto ${\cal D}$}

   \section{\bf When  $X$  is symmetric and $\int_0^1{1\over H^2(t)} dt<+\infty $ } 
  
 From now on, we suppose that $X$ is symmetric. 
Let 
$$s(x): =\int_0^x{1\over H^2(t)} dt $$
 and $\psi$ be the exponent of the L\'evy process :
  $${\bf P}(e^{-iuX_t})=: e^{-t\psi(u)}$$
 \medskip
  {\bf Theorem} 
{\it   We suppose that   $s(x)<+\infty$ for every $x>0$.

  1) There exists a unique string $m\in {\cal M}$ such that  
   $$m(s(x)) =\int_0^x{H^4(t)\over 4}dt$$
   
2)  The $A$-solution associated to $m$ is given by :

$$A(s(x),\lambda)={e^{\lambda {x\over 2}}\phi(x,\lambda )+e^{-\lambda {x\over 2}}{\phi}(x,-\lambda)\over H(x)}$$

3) The spectral  measure associated to $m$ is $(\psi(u)+1)du$ in  the first case and $\psi(u)du$ in  the second case.
}
\medskip

{\bf Argument} : 1) is trivial.

 2) is a consequence of elementary transformations of the differential system given in proposition 3. 
 
 In next section, we give a summary  of the proof of 3) in the first case. The second case follows when varying the rate of killing (replace  the rate $1$ of the independant exponential time when the L\'evy process $X$ is killed by $\mu$, any positive constant) and then make $\mu$ go to $0$.
\medskip

{\bf Proposition}
{\it The following four  assertions are equivalent :

$$(i)\qquad \psi(u)du \in {\cal D}$$
$$(i') \qquad (\psi(u)+1) du \in {\cal D}$$
$$(ii) \qquad \hbox{ the paths of $X$ have bounded variation}$$ 
$$(iii)\qquad  \int_0^1 {1\over H^2(t)}dt<+\infty$$
}

{\bf Argument}  The equivalence between the properties $(i)$ and $(i')$ is trivial.
the equivalence between $(i)$ and $(ii)$ has been proved by Vigon in his thesis.

The equivalence between $(i)$ and $(iii)$ follows easily from the previous theorem.

\section {\bf Proof of part 3 of the theorem }

{\bf Lemma 1 } {\it In the first case and for any L\'evy process, we have for $\lambda >0$
$$-\log  {\bf Q}^{\uparrow}(e^{-\lambda .F})={1\over 2\pi} \int_{-\infty}^{+\infty}\log(\psi(u)+1)({1\over -iu+\lambda }-{1\over  -iu})du$$
and if $X$ is symmetric :
$$-\log{\bf Q}^{\uparrow}(e^{-\lambda .F})={2\lambda \over \pi} \int_0^{+\infty} {\log(\psi(u)+1)\over u^2+\lambda ^2}du$$
}
\medskip
 
{\bf Argument}  for $\lambda \in i R$, the independence property  of trajectories before and after the minimum (see  first section) gives  the equality :  
$$-\log  {\bf Q}(e^{-\lambda  F})=-\log{\bf Q}^{\uparrow}(e^{-\lambda  F})-\log \check {\bf Q}^{\uparrow}(e^{\lambda  F})$$

One can prove that $\lambda \mapsto -\log {\bf Q}(e^{-\lambda  F})$ is the  the exponent of a L\'evy process with bounded variation and no drift and that  the above equality is its decomposition into a sum of the exponent of a subordinator and the exponent of the opposite of a subordinator.
Then extract  $-\log{\bf Q}^{\uparrow}(e^{-\lambda  F})$ from $-\log {\bf Q}(e^{-\lambda  F})$ by a (sort of) Stieljes transform.
\medskip

{ \bf Lemma 2} {\it [Krein's Entropy Formula] When $\Delta(du)=\Delta'(u)du$, one gets  
$${2\lambda \over \pi}\int_0^{+\infty}{\log[\Delta'(u)]\over u^2+\lambda ^2}du =\lim_{x\to +\infty}\log(-e^{2\lambda .s^{-1}(x)}D'(x,\lambda )D(x,\lambda))+\log \lambda $$
}
{\bf Final Argument}   
We have noticed previously that 
$$D(x)=A(x)\int_x^{+\infty} {1\over A^2(t)}dt \quad \hbox{ and }\quad  A(s(x))\sim e^{\lambda  {x\over 2}}{\bf Q}^{\uparrow}(e^{-\lambda  F}) \quad (x\to +\infty )$$

With the help of these two properties, one can easily find out that the second member of lemma 2 is equal to $-\log {\bf Q}^{\uparrow}(e^{-\lambda  F})$. Lemma 1  then gives the equality $\Delta'(u)= \psi(u)+1$  by identification.

\section{ When $X$ is symmetric and has unbounded variations}
In that case, we have   
$s(x)=\int_0^x{1\over H^2(t)}dt=+\infty $.
But one  still  can define : 
$$\tilde A(x):={e^{\lambda {x\over 2}}\phi(x,\lambda )+e^{-\lambda {x\over 2}}\phi (x,-\lambda )\over H(x)}$$
$${\tilde D}(x):={\tilde A}(x)\int_x^{+\infty} {1\over {\tilde A}^2(t)} {1\over H^2(t)} dt$$
Notice that ${\tilde D}(x)\in ]0,+\infty[$ for every $x>0$ but ${\tilde D}(0)=+\infty$.

Still  we have that the measure ${\psi(u)+1\over u^2+1}du\hbox{ is in }{\cal D}$
 (this  is true for any symmetric L\'evy process ).
Denote  by $D_1$  the $D$-solution associated to this  spectral measure ${\psi(u)+1\over u^2+1}du$.

And  define $$t(x):=x+{1\over H^2(x)}.{{\tilde D}(x,1)\over - {\tilde D}'(x,1)}$$

{\bf Theorem} {\it [ Adaptation of `` Rule 4" of Dym-Mac Kean ] : 

$t$ is increasing,\  $t(0^+)=0$\,  and $t(+\infty )=+\infty$ and, 

$${D_1(t(x),\lambda )\over  -{\tilde D}'(x,\lambda )}={1\over \lambda ^2-1} \biggl ({{\tilde D}(x,\lambda )\over -{\tilde D}'(x,\lambda)}-{{\tilde D}(x,1)\over -{\tilde D}'(x,1)}\biggr)$$
}

 {\bf  Remark} 
 
 The Brownian case can be computed through  the solution corresponding to the spectral measure  $\Delta(du)=du \in {\cal D}$!!!
 \bigskip

 \centerline {\bf Bibliography}
 
{\bf Bertoin} (1996) {\it L\'evy Processes }, Cambridge Univ.Press, Cambridge. 

{\bf Dym-Mac Kean} (1976) {\it Gaussian Processes, Function Theory, and the Inverse Spectral Problem},
Academic Press

{\bf Rogozin} (1972) : {\it  The distribution of the first hit for stable  and asymptotic stable walks on interval.} Theor. 
Probab. Appl. 17, 332-338
 
{\bf Sato} (1999) {\it L\'evy Processes and Infinitely Divisible Distributions}, Cambridge studies in advanced mathematics

{\bf Vigon} (2002) { \it Simplifiez vos L\'evy en titillant la factorisation de Wiener-Hopf}. Thesis.

\end{document}